\documentclass[10pt]{amsart}%
\usepackage{amssymb}
\usepackage{amsmath,amsthm}
\usepackage{graphicx}
\usepackage{hyperref}
\usepackage{url}
\usepackage{amsfonts}
\usepackage{multicol}
\usepackage{amsmath}
\usepackage{verbatim}%
\providecommand{\U}[1]{\protect\rule{.1in}{.1in}}
\theoremstyle{theorem}
\newtheorem{theorem}{Theorem}
\theoremstyle{definition}

\newtheorem{exercise}{Exercise}
\allowdisplaybreaks
\makeatletter
\@addtoreset{footnote}{page}
\makeatother
\begin{document}
\title{Generalized Sine Functions, Complexified}
\author{Pisheng Ding}
\address{Department of Mathematics \\
Illinois State University \\
Normal, IL 61790, USA}
\email{pding@ilstu.edu}
\author{Sunil K. Chebolu}
\address{Department of Mathematics \\
Illinois State University \\
Normal, IL 61790, USA}
\email{schebol@ilstu.edu}
\thanks{The second author is supported by Simons Foundation: Collaboration Grant for
Mathematicians (516354). }

\begin{abstract}
Generalized sine and cosine functions, $\sin_{n}$ and $\cos_{n}$, that
parametrize the generalized unit circle $x^{n}+y^{n}=1$ are, much like their
classical circular counterparts, extendable as complex analytic functions. In
this article, we identify the natural domain on which $\sin_{n}$ is a
conformal equivalence from a polygon to the complex plane with $n$ slits. We
also give some geometric and analytic applications.

\end{abstract}
\maketitle

\section{Introduction}

A 2011 article \cite{Wood} gives an account of certain generalized
trigonometric functions that are suited for describing generalized circles
such as the curve $x^{4}+y^{4}=1$. Its author hopes that his piece
\textquotedblleft\lbrack offer] many opportunities for the reader to extend
the theory into studies of $\cdots$ complex analysis\textquotedblright.
Indeed, extending this story into the imaginary realm is our express
objective, and we also hope to offer opportunities for the reader to gain
further insight into the complex (pun intended)\ behavior of these real functions.

These generalized trigonometric functions seem to be first systematically
introduced by a 1959 paper \cite{Shelupsky} and interest in them has since
grown. As recently as 2021, the second author investigated aspects of them
with his students in an undergraduate research course; see \cite{Cheb}. The
study of these functions, dubbed (by \cite{Wood} and \cite{Poodiack})
\textquotedblleft squigonometry\textquotedblright\ in recreational
mathematics, is now an active field of pure and applied analysis due in part
to its connection to a family of differential operators known as the
$(p,q)$-Laplacian; see, e.g., \cite{Lang-Edmunds}. On the most basic level,
these functions resemble, but also differ significantly from, their circular
counterparts. For background, we give a brief introduction to their
definitions and immediate properties.

We first recall the familiar fact that, for $y\in(-1,1)$,%
\[
\arcsin y=\int_{0}^{y}\frac{1}{\sqrt{1-x^{2}}}dx\,\text{.}%
\]
In other words, the function $\sin:(-\pi/2,\pi/2)\rightarrow(-1,1)$ is the
inverse of the function $y\mapsto\int_{0}^{y}(1-x^{2})^{-1/2}dx$. As is
well-known, $\sin t$ is the $y$-coordinate of a point on the unit circle with
$t$ signifying either the length of an arc or twice the area of a circular sector.

Now, let $n$ be an integer with $n>2$. For $y\in(a_{n},1)$ where
$a_{n}=-\infty$ for $n$ odd and $a_{n}=-1$ for $n$ even, consider%
\[
F_{n}(y):=\int_{0}^{y}\frac{1}{(1-x^{n})^{(n-1)/n}}dx\text{\thinspace,}%
\]
which is strictly increasing and hence invertible. Cued by the fact that
$\arcsin1=\pi/2$, we define $\pi_{n}$ so that $F_{n}(1)=\pi_{n}/2$. Then,
define $\sin_{n}$ to be $F_{n}^{-1}:(F_{n}(a_{n}),\pi_{n}/2)\rightarrow
(a_{n},1)$ and $\cos_{n}t$ to be $\sqrt[n]{1-\sin_{n}^{n}t}$. By definition,%
\begin{equation}
\sin_{n}^{n}t+\cos_{n}^{n}t=1\text{\thinspace.} \label{Pyth}%
\end{equation}
Thus, $t\mapsto(\cos_{n}t,\sin_{n}t)$ parameterizes the right-half-plane
portion of the curve%
\[
C_{n}:=\{(x,y)\mid x^{n}+y^{n}=1\}\,\text{.}%
\]

\begin{figure}[h]
\centering\includegraphics[scale=0.4]{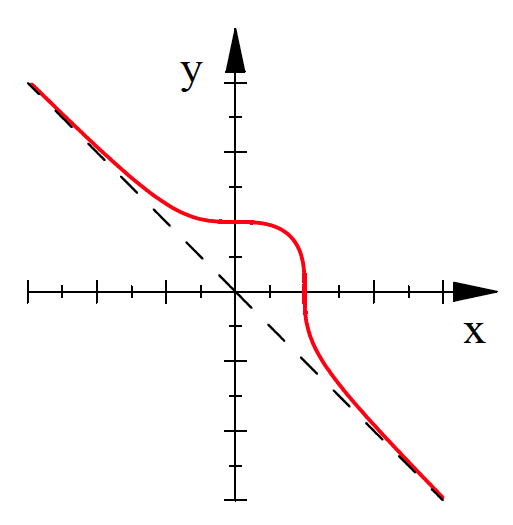}\caption{{}$C_{3}$ and its
asymptote.}%
\label{C3}%
\end{figure}

\begin{figure}[h]
\centering\includegraphics[scale=0.4]{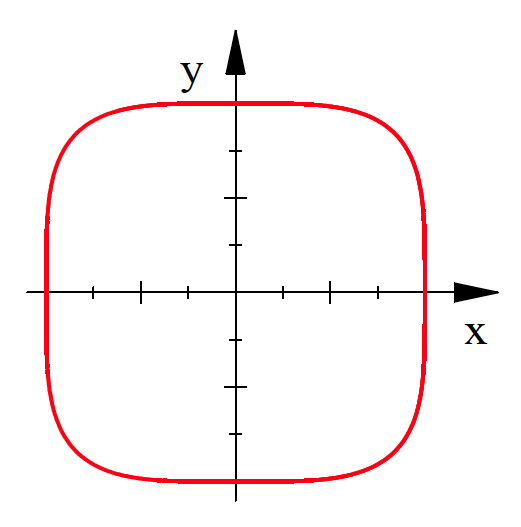}\caption{{}$C_{4}$.}%
\label{C4}%
\end{figure}

We seek to interpret the parameter $t$. It can be shown from the definition
that%
\begin{equation}
\sin_{n}^{\prime}t=\cos_{n}^{n-1}t\text{\quad and\quad}\cos_{n}^{\prime
}t=-\sin_{n}^{n-1}t\text{\thinspace.} \label{ODE}%
\end{equation}
Consider now a sector bounded by two rays from the origin and an infinitesimal
arc on $C_{n}$ due to an infinitesimal increment $dt$ in the parameter $t$. By
(\ref{Pyth}) and (\ref{ODE}), this infinitesimal sector has differential area%
\begin{align*}
\frac{1}{2}r^{2}d\theta &  =\frac{1}{2}(x^{2}+y^{2})\cdot d\left(
\arctan\frac{y}{x}\right)  =\frac{1}{2}\left(  xdy-ydx\right) \\
&  =\frac{1}{2}\left[  \cos_{n}t\cdot\cos_{n}^{n-1}t-\sin_{n}t\cdot(-\sin
_{n}^{n-1}t)\right]  dt=\frac{1}{2}dt\,\text{.}%
\end{align*}
This reveals the geometric significance of $t$:\ if $A$ denotes the point
$(1,0)$ and $P(t)$ the point $(\cos_{n}t,\sin_{n}t)$, then $t$ is
\textit{twice} the (signed) area of the sector $OAP(t)$; as a result, $\pi
_{n}/4$ is then the area in the first quadrant bounded by $C_{n}$. This
geometric interpretation of $t$ allows the extension of the domain of
$\cos_{n}$ and $\sin_{n}$ beyond the interval $(F_{n}(a_{n}),\pi_{n}/2]$ as
follows. With $P(0)=A$, let $P(t)$ describe $C_{n}$ with increasing angular
coordinate as $t$ increases such that the sector $OAP(t)$ has\ area $t/2$. The
pair $(\cos_{n}t,\sin_{n}t)$ are then defined to be the Cartesian coordinates
of $P(t)$. Thus, the area interpretation for $t$ in $\sin t$ survives the
generalization. (In the case $n=4$, \cite{Levin} gives an interesting
application of this area interpretation.) Figures \ref{C3} and \ref{C4} show
$C_{3}$ and $C_{4}$, from which it is evident that $\sin_{n}$ is odd iff $n$
is even.



Having made sense of $\sin_{n}t$ and $\cos_{n}t$ for (certain) real $t$, how
on earth can we attach meaning to $\sin_{n}z$ and $\cos_{n}z$ for a complex
input $z$? There is a simple approach, but it will not get us very far. To
appreciate its limitation, let us briefly describe it. As $F_{n}$ can be
expanded into (and is equal to) a convergent power series at $0$ with
$F^{\prime}(0)\neq0$, its inverse $\sin_{n}$ is also analytic at $0$.
Repeatedly differentiating the identity $\sin_{n}(F_{n}(y))=y$, we can find
(in principle)\ $\sin_{n}^{(k)}(0)$ for $k\geq1$. Writing down the Maclaurin
series for $\sin_{n}$, we can force it to accept a complex input $z$ that is
within a certain distance (the radius of convergence)\ from $0$. For example,
we find in this way%
\[
\sin_{3}z=z-\frac{1}{6}z^{4}+\frac{2}{63}z^{7}-\frac{13}{2268}z^{10}+\frac
{23}{22113}z^{13}-\cdots\text{ .}%
\]
(For a discussion of the Maclaurin series for $\sin_{n}$, see, e.g.,
\cite{Cheb}.) But this is not very illuminating concerning the behavior of
$\sin_{3}$ as a complex function. With this approach, we do not even know its
radius of convergence!

We aim to identify $\sin_{n}$'s natural domain of analyticity and to describe
$\sin_{n}$ as a conformal equivalence. We first treat $\sin_{4}$; this special
case will indicate how $\sin_{n}$ behaves in general. As we shall see, our
theory will have tangible \textit{real} (pun intended)\ implications. For a
preview, we assert a claim that we will effortlessly prove:\ \textit{The area
between }$C_{3}$\textit{\ and its asymptote is trisected by the coordinate
axes}.

The exposition assumes some basic notions in complex analysis that can be
found in any undergraduate textbook on the subject, e.g., \cite{Bak}. For
students learning the subject, this article provides enriching applications of
many key ideas.

\section{A special case:\ sin$_{4}$}

We introduce notation that will be convenient for exposition throughout the article.

A point in the plane is either expressed as a complex number or named by an
uppercase letter. Given two points $X$ and $Y$, we write $[X,Y]$ for the
(closed)\ line segment between them and $[X,Y)$ for the half-open segment that
omits $Y$. In either case, we write $|XY|$ for the length of the segment. We
let $I=[0,1]$ and $J=[1,\infty)$.

For $a\in%
\mathbb{C}
$ and $S\subset%
\mathbb{C}
$, we write $aS$ to denote $\{az\mid z\in S\}$.

Fix $n\geq3$. Let $\omega_{n}=e^{i2\pi/n}$. Let $A_{n}=\pi_{n}/2$ and
$B_{n}=\omega_{n}A_{n}$. Define $V_{n}$ to be the wedge-shaped region
$\{re^{i\theta}\mid r>0;\,\theta\in(0,2\pi/n)\}$.

\begin{figure}[h]
\centering\includegraphics[scale=0.2]{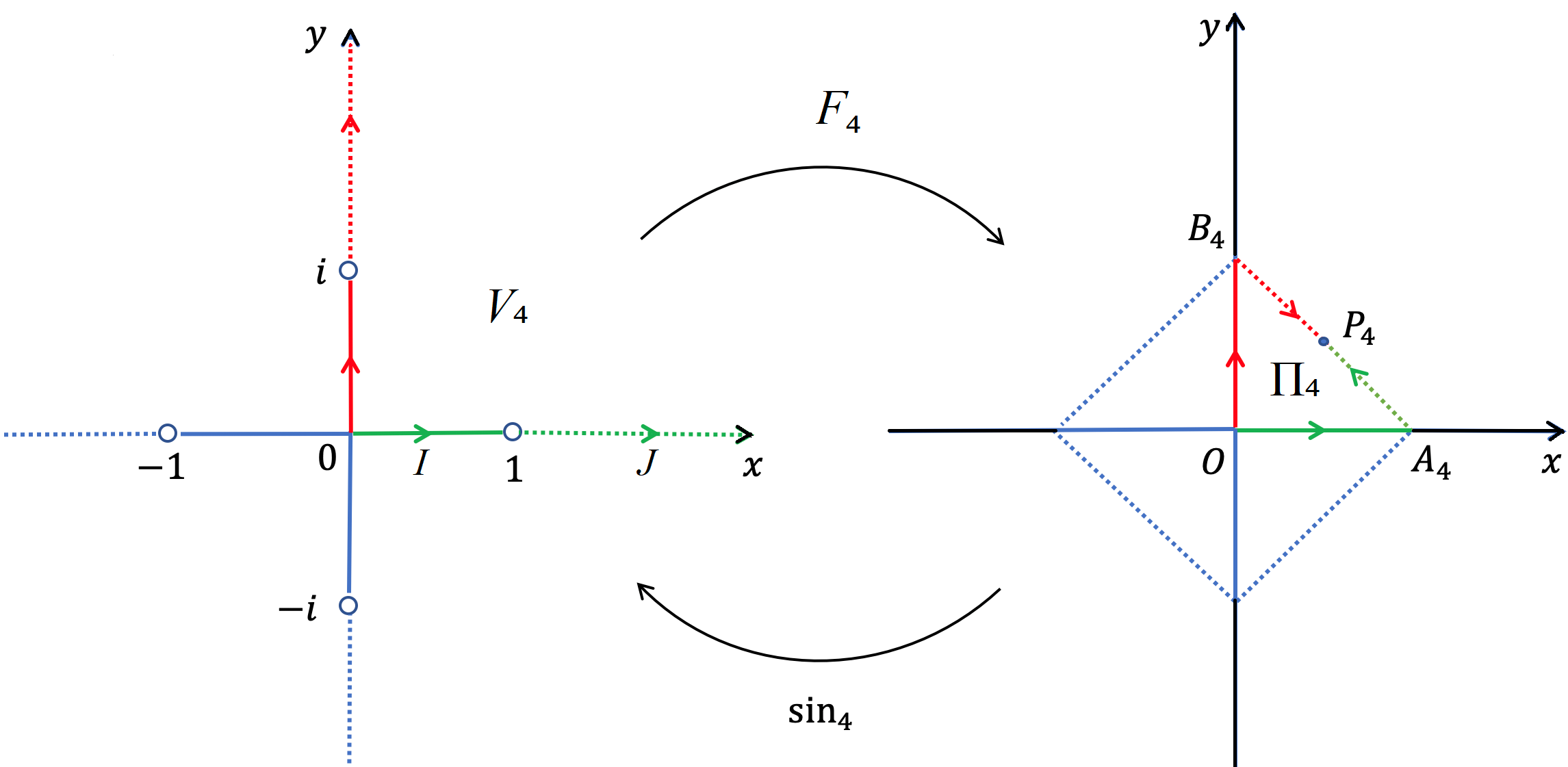}\caption{{}$F_{4}$ and
$\sin_{4}$.}%
\label{n=4}%
\end{figure}

We now turn to the case $n=4$, in which $\omega_{4}=i$ and $V_{4}$ is the open
first quadrant. Consider, for $z\in V_{4}$,
\[
K_{4}(z)=\frac{1}{(1-z^{4})^{3/4}}%
\]
with the requirements that $K_{4}(0)=1$ (the principal branch of the power)
and that $K_{4}$ be continuous on $V_{4}$. As $V_{4}$ is free of singularities
of $K_{4}$, $K_{4}$ is analytic on $V_{4}$; as $V_{4}$ is simply connected,
$K_{4}$ has a primitive defined by a path-independent integral%
\[
F_{4}(z)=\int_{0}^{z}K_{4}(\zeta)d\zeta\text{\thinspace.}%
\]
We first analyze the behavior of $F_{4}$ on $\partial V_{4}$. Refer to
Figure\thinspace\ref{n=4} as you follow the analysis.

\begin{itemize}
\item For $x\in\lbrack0,1)$,%
\[
K_{4}(x)=\frac{1}{\sqrt[4]{1-x^{4}}^{3}}\text{\quad and\quad}F_{4}(x)=\int%
_{0}^{x}\frac{1}{\sqrt[4]{1-t^{4}}^{3}}dt
\]
where $\sqrt[4]{\quad}$ is the usual principal $4$th root of a
\textit{positive} number. While $K_{4}(1)$ is \textit{not} defined, $F_{4}(1)$
\textit{is}, because it equals a convergent improper integral. Thus,
$F_{4}[I]=[0,F_{4}(1)]=[O,A_{4}]$.

\item For $x>1$,%
\[
K_{4}(x)=e^{i\pi3/4}\frac{1}{\sqrt[4]{x^{4}-1}^{3}}%
\]
where the phase factor $e^{i\pi3/4}$ is due to the requirement that $1/K_{4}$
(the denominator of $K_{4}$) be continuous in the closure $\overline{V}_{4}$
of $V_{4}$\footnote{For $z$ near $1$, factor $(1-z^{4})^{3/4}$ into
$g(z)(z-1)^{3/4}$ where $g$ is continuous and nonzero at $1$. Thus, it
suffices to track the phase change of $(z-1)^{3/4}$ as $z$ traverses an arc in
$\overline{V}_{4}$ from $1-\epsilon$ to $1+\epsilon$.}. Then,%
\[
F_{4}(x)=F_{4}(1)+\int_{1}^{x}K_{4}(t)dt=\frac{\pi_{4}}{2}+e^{i\pi3/4}\int%
_{1}^{x}\frac{1}{\sqrt[4]{t^{4}-1}^{3}}dt\text{ ,}%
\]
showing that $F_{4}[J]$ is a segment $[A_{4},P_{4})$ with an angle of
inclination $3\pi/4$ where $P_{4}=\lim_{x\rightarrow\infty}F(x)$. Note that
$|A_{4}P_{4}|=\int_{1}^{\infty}1/\sqrt[4]{t^{4}-1}^{3}dt$, a \textquotedblleft
doubly\textquotedblright\ improper integral that is convergent.

\item For $z=it$ with $t\in\lbrack0,1)$,%
\[
K_{4}(z)=(1-(it)^{4})^{-3/4}=\left(  \sqrt[4]{1-t^{4}}\right)  ^{-3}%
\]
and%
\[
F_{4}(z)=\int_{0}^{t}K_{4}(i\tau)id\tau=i\int_{0}^{t}\frac{1}{\sqrt[4]%
{1-\tau^{4}}^{3}}d\tau\text{ ,}%
\]
which shows that $F_{4}[iI]=[O,B_{4}]$.

\item For $z=it$ with $t>1$,%
\[
K_{4}(z)=e^{-i\pi3/4}\left(  \sqrt[4]{t^{4}-1}\right)  ^{-3}\text{ .}%
\]
The phase factor $e^{-i\pi4/3}$ is again dictated by the continuity of
$1/K_{4}$ in\textit{\ }$\overline{V}_{4}$. Thus,%
\[
F_{4}(z)=F_{4}(i)+\int_{1}^{t}K_{4}(i\tau)id\tau=\frac{\pi_{n}}{2}%
i+e^{-i\pi3/4}i\int_{1}^{t}\frac{1}{\sqrt[4]{\tau^{4}-1}^{3}}d\tau
\text{\thinspace,}%
\]
showing that $F_{n}[iJ]$ is a segment $[B_{4},P_{4}^{\prime})$ that makes an
angle of \textit{signed} measure $-3\pi/4$ measured from (the positive
direction of) the $y$-axis. While we shall next argue that $P_{4}^{\prime
}=P_{4}$, note that $A_{4},P_{4},P_{4}^{\prime},B_{4}$ are colinear and that%
\begin{equation}
|B_{4}P_{4}^{\prime}|=\int_{1}^{\infty}\frac{1}{\sqrt[4]{t^{4}-1}^{3}%
}dt=|A_{4}P_{4}|\text{ .}\label{Integral}%
\end{equation}

\item We make a simple but important claim, which will imply that
$P_{4}^{\prime}$ coincides with $P_{4}$ and hence is the midpoint of
$[A_{4},B_{4}]$:%
\begin{equation}
\lim_{|z|\rightarrow\infty}F_{4}(z)=P_{4}\text{\thinspace,}\label{F(infinity)}%
\end{equation}
where the limit is taken with $z\in\overline{V}_{4}$. To verify
(\ref{F(infinity)}), note that $F_{4}(re^{i\theta})-F_{4}(r)$ is expressible
as the integral of $F_{4}^{\prime}=K_{4}$ along a circular arc. In detail, let
$\epsilon>0$ be given. We leave to the reader to argue that the modulus of%
\[
F_{4}(re^{i\theta})-F_{4}(r)=\int_{z=re^{it},\;t\in\left[  0,\theta\right]
}\frac{1}{(1-z^{4})^{3/4}}dz
\]
can be bounded by $\epsilon/2$ for sufficently large $r$ and all $\theta
\in\left[  0,2\pi/n\right]  $. Since $P_{4}:=\lim_{r\rightarrow\infty}%
F_{4}(r)$, we have $|F_{4}(r)-P_{4}|<\epsilon/2$ for sufficiently large $r$.
Then,%
\[
|F_{4}(re^{i\theta})-P_{4}|<\epsilon
\]
for sufficently large $r$ and all $\theta\in\left[  0,2\pi/n\right]  $,
establishing (\ref{F(infinity)}).
\end{itemize}

Letting $\Pi_{4}$ denote the closed isosceles right triangle $\triangle
OA_{4}B_{4}$, we summarize our key findings concerning $F_{4}$ on $\partial
V_{4}$:

\begin{quote}
$F_{4}$\textit{\ maps }$\partial V_{4}$\textit{\ homeomorphically onto
}$\partial\Pi_{4}\setminus\{P_{4}\}$\textit{\ with }$F_{4}(0)=0,$%
\textit{\ }$F_{4}(1)=A_{4},$\textit{\ }$F_{4}(i)=B_{4}$\textit{, while }%
$\lim_{z\in\overline{V}_{4};\text{\thinspace}|z|\rightarrow\infty}%
F_{4}(z)=P_{4}$.
\end{quote}

Turning to $F_{4}$ on $V_{4}$, we argue, by the argument principle, that
$F_{4}$ maps $V_{4}$ one-to-one onto $\mathring{\Pi}_{4}$ (the interior of
$\Pi_{4}$). Let $w\in\mathring{\Pi}_{4}$. For $R>0$, let $\Gamma_{R}$ be the
positively oriented boundary of the circular sector with vertices $0$, $R $,
and $iR$. The preceding analysis implies that, for all sufficiently large $R$,
$F_{4}[\Gamma_{R}]$ winds around $w$ exactly once, proving that $w$ is
attained as a value of $F_{4}$ exactly once (with multiplicity 1). Thus, being
a bijective analytic map, $F_{4}:V_{4}\rightarrow\mathring{\Pi}_{4} $ is a
conformal equivalence.

At long last, we are able to define $\sin_{4}:\mathring{\Pi}_{4}\rightarrow
V_{4}$ as the inverse of $F_{4}$. As $\sin_{4}$ (or, to be precise, its
boundary extension) maps $[O,B_{4})$ into the segment $[0,i)$, we can apply
Schwarz reflection principle to define $\sin_{4}$ analytically on
$i\mathring{\Pi}_{4}$, the reflected image of $\mathring{\Pi}_{4}$ across
$[O,B_{4})$.\footnote{Schwarz reflection principle also guarantees that the
extended function is analytic on $[O,B_{4})$.} Repeating the same argument two
more times extends $\sin_{4}$ analytically on the interior $\Omega_{4}$ of the
union of the four (closed) triangles $i^{k}\Pi_{4}$, $k\in\{0,1,2,3\}$. The
extended function $\sin_{4}$ maps each (open)\ triangle $i^{k}\mathring{\Pi
}_{4}$ conformally onto the quadrant $i^{k}V_{4}$; its boundary extension maps
each vertex $i^{k}A_{4}$ to $i^{k}$, while it maps, in a two-to-one fashion,
$\overline{\Omega}_{4}$'s four sides with the midpoints and vertices removed
onto the four (open) rays $\cup_{k=0}^{3}i^{k}\mathring{J}$.

We summarize the properties of $\sin_{4}$ as follows, where $\Sigma_{4}$
denotes $%
\mathbb{C}
\setminus\left(  \cup_{k=0}^{3}i^{k}J\right)  $, the plane with four slits.

\begin{enumerate}
\item $\sin_{4}:\Omega_{4}\rightarrow\Sigma_{4}$\textit{\ is a bijective
analytic map and hence a conformal equivalence.}

\item \textit{For }$k\in\{0,1,2,3\}$\textit{, }$|\sin_{4}z|\rightarrow\infty$
as $z\rightarrow i^{k}P_{4}$ in $\Omega_{4}$.

\item \textit{The continuous extension of }$\sin_{4}$\textit{\ maps }%
$i^{k}A_{4}$\textit{ to }$i^{k}$\textit{ and maps }$\partial\Omega
_{4}\setminus\{i^{k}P_{4},i^{k}A_{4}\}_{k=0}^{3}$\textit{\ two-to-one onto
}$\cup_{k=0}^{3}i^{k}\mathring{J}$\textit{.}
\end{enumerate}

With this background, the general case can be conquered with ease.

\section{sin$_{n}$ for $n\geq3$}

The notations introduced in the preceding section remain in force.

Much like the case $n=4$, we consider, for $z\in V_{n}$,%
\[
K_{n}(z)=\frac{1}{(1-z^{n})^{(n-1)/n}}\text{\quad and\quad}F_{n}(z)=\int%
_{0}^{z}K_{n}(\zeta)d\zeta\text{.}%
\]

\begin{figure}[h]
\centering\includegraphics[scale=0.2]{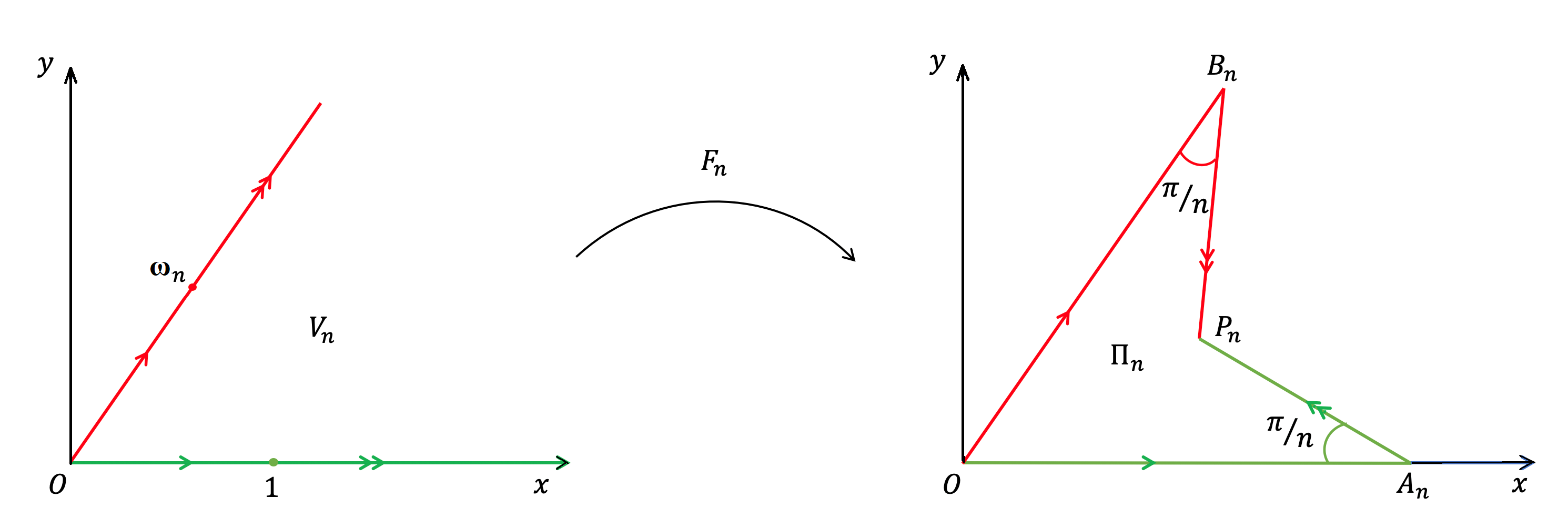}\caption{{}$F_{n}%
:V_{n}\rightarrow\Pi_{n}$.}%
\label{Fn}%
\end{figure}

\begin{figure}[h]
\centering\includegraphics[scale=0.2]{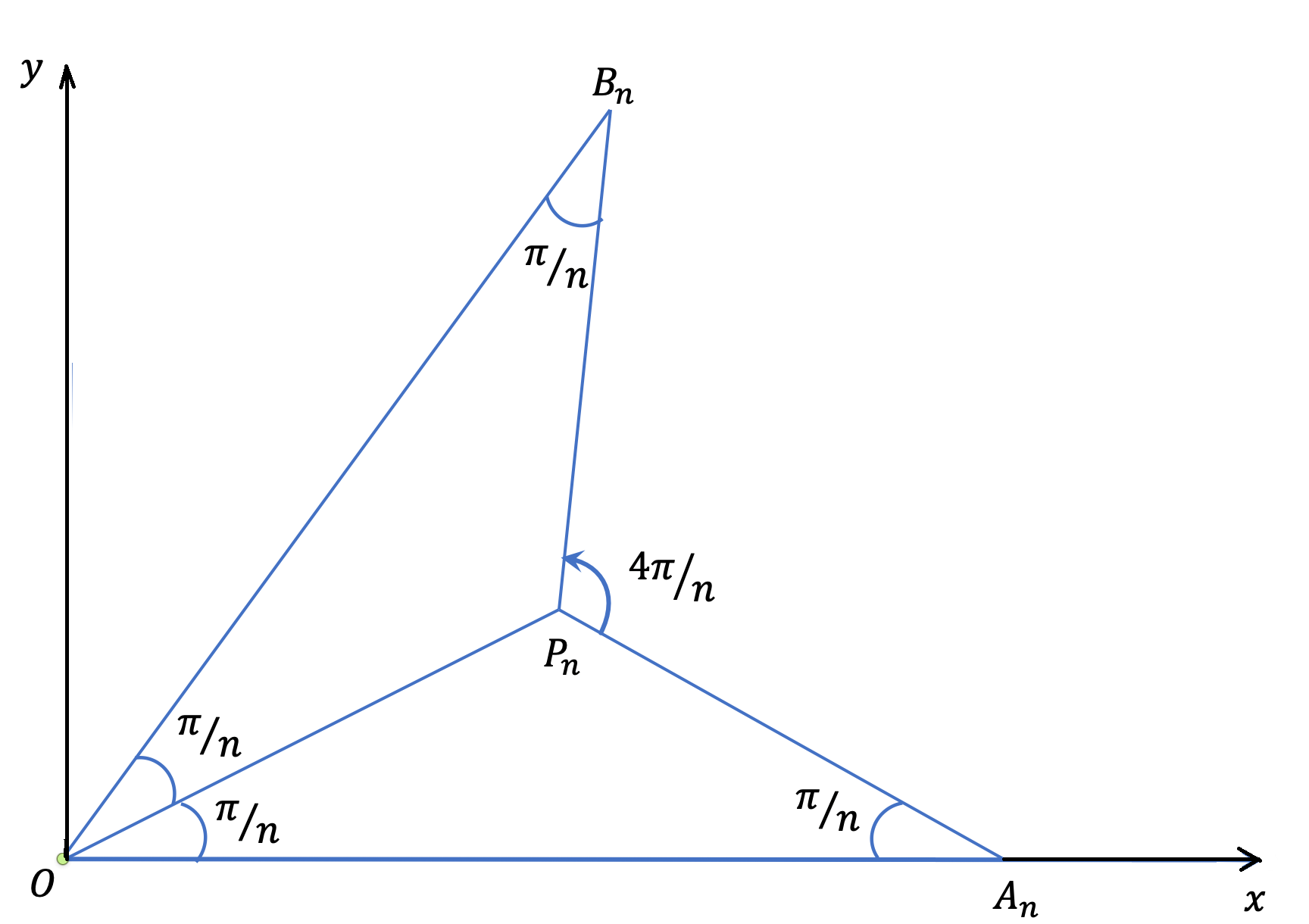}\caption{{}The polygon
$\Pi_{n}$.}%
\label{PI-n}%
\end{figure}

Calculation similar to that in the case $n=4$ reveals the action of $F_{n}$ on
$\overline{V}_{n}$, as depicted in Figure\thinspace\ref{Fn}. While leaving the
details to the reader, we note a difference: when $n\neq4$, the three points
$A_{n},P_{n},B_{n}$ are no longer colinear, and $\triangle OA_{4}B_{4}$
becomes a (not necessarily convex)\ quadrilateral $\Pi_{n}$ with successive
vertices $O,A_{n},P_{n},B_{n}$. In fact, $\Pi_{n}$ admits a unified
description for all $n\geq3$, which we now give; see Figure\thinspace
\ref{PI-n}.

Let $P_{n}$ denote the point on the bisector of $\angle A_{n}OB_{n}$ such that
$\measuredangle OA_{n}P_{n}=\pi/n$. Then $\Pi_{n}$ is the closed polygon
$OA_{n}P_{n}B_{n}$ (with the caveat that, in $\Pi_{4}$, $P_{4}$ becomes a
degenerate vertex.) Note that the \textit{counterclockwise} angle from
$\overrightarrow{P_{n}A_{n}}$ to $\overrightarrow{P_{n}B_{n}}$ has measure
$4\pi/n$, which is less than $\pi$ iff $n>4$. Consequently, $\Pi_{n}$ is
nonconvex iff $n>4$, whereas $\Pi_{3}$ is an equilateral parallelogram and
$\Pi_{4}$ is an isosceles right triangle (already seen). Simple trigonometry
gives the location of $P_{n}$:%
\begin{equation}
P_{n}=\left(  \frac{\pi_{n}}{4}\sec\frac{\pi}{n}\right)  e^{i\pi/n}\text{ .}
\label{Pn}%
\end{equation}

Let $\Omega_{n}$ denote the interior of $\cup_{k=0}^{n-1}\omega_{n}^{k}\Pi
_{n}$ (the union of the images of $\Pi_{n}$ under rotations by angle $2\pi
k/n$ for various $k$). Note that $\Omega_{n}$ is a (open)\ equilateral
$2n$-gon except that $\Omega_{4}$ is a (open) square. See Figure\thinspace
\ref{Domains(3,4,6)} for $\Omega_{3}$, $\Omega_{4}$, and $\Omega_{6}$.

\begin{figure}[h]
\centering\includegraphics[scale=0.2]{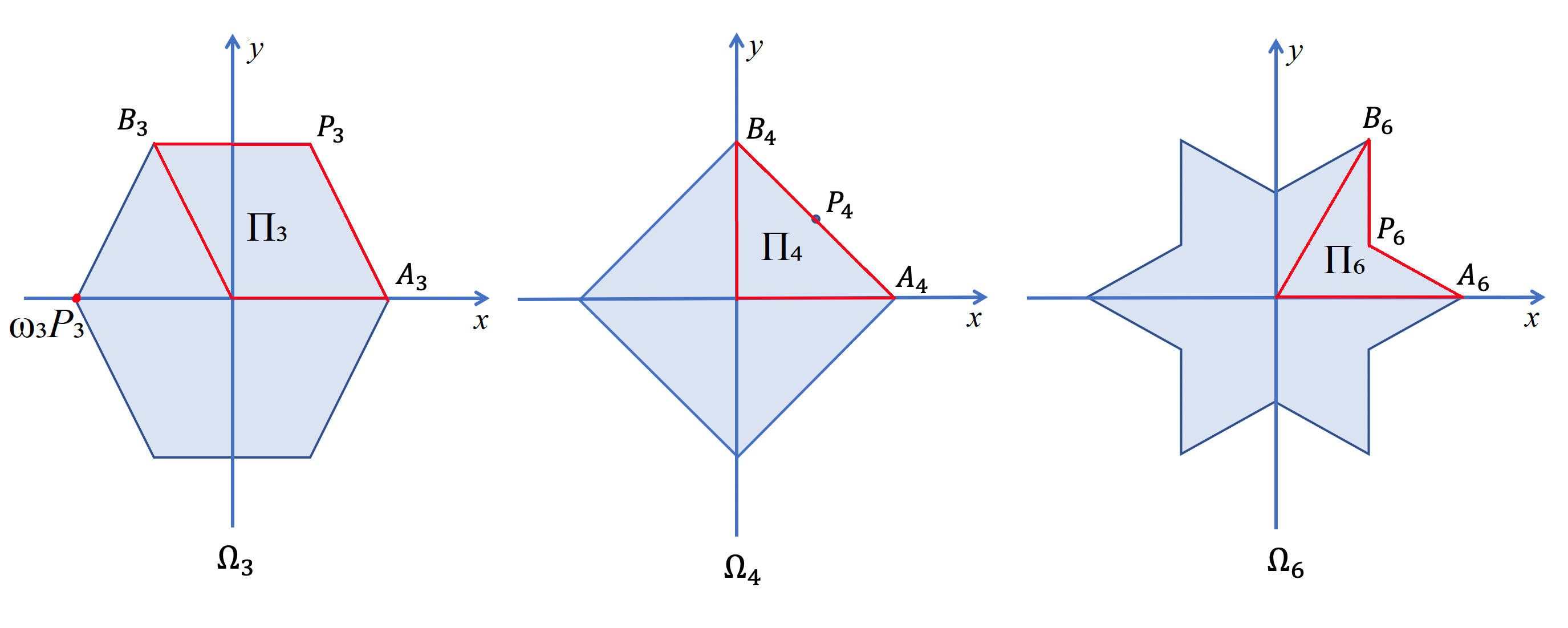}\caption{Domains for {}%
$\sin_{3}$, $\sin_{4}$, and $\sin_{6}$.}%
\label{Domains(3,4,6)}%
\end{figure}

Lastly, let $\Sigma_{n}$ denote $%
\mathbb{C}
\setminus\left(  \cup_{k=0}^{n-1}\omega_{n}^{k}J\right)  $, the plane with $n$ slits.

Parallel to the case $n=4$, we have our main result.

\begin{theorem}
Let $n\geq3$ and $k\in\{0,1,\cdots,n-1\}$.

\begin{enumerate}
\item $\sin_{n}:\Omega_{n}\rightarrow\Sigma_{n}$ is a bijective analytic map
and hence a conformal equivalence.

\item $|\sin_{n}z|\rightarrow\infty$ \textit{as }$z\rightarrow\omega_{n}%
^{k}P_{n}$ in $\Omega_{n}$.

\item $\sin_{n}$, continuously extended, maps $\omega_{n}^{k}A_{n}$ to
$\omega_{n}^{k}$ and maps $\partial\Omega_{n}\setminus\{\omega_{n}^{k}%
A_{n},\omega_{n}^{k}P_{n}\}_{k=0}^{n-1}$ two-to-one onto $\cup_{k=0}%
^{n-1}\omega_{n}^{k}\mathring{J}$.
\end{enumerate}
\end{theorem}

\section{Real consequences}

As promised, we explore some \textit{real} applications of our theory.

\subsection{Radius of convergence}

Earlier, we wondered how we might find the radius of convergence $R_{n}$ of
$\sin_{n}$'s Maclaurin series. Having extended the domain of $\sin_{n}$, we
now have an exact answer. Recall that $D(0;R_{n})$ (the origin-centered open
disc of radius $R_{n}$) is the largest among all origin-centered discs on
which $\sin_{n}$ is analytic. Except for the case $n=3$, it should be clear
from Figures \ref{PI-n} and \ref{Domains(3,4,6)} that the obstruction to
analyticity closest to the origin occurs at $P_{n}$ (and at $\omega_{n}%
^{k}P_{n}$). Thus, for $n>3$, $R_{n}=|OP_{n}|$; by (\ref{Pn}), we have%
\begin{equation}
R_{n}=\frac{\pi_{n}}{4}\sec\frac{\pi}{n}\text{\thinspace.} \label{Rn}%
\end{equation}

Let us turn to the case $n=3$. As $D(0;|OP_{3}|)\nsubseteq\Omega_{3}$, the
question is then whether $\sin_{3}$ can be analytically defined on
$D(0;|OP_{3}|)$.\ Answer:\ Yes, and consequently, formula (\ref{Rn}) holds
without exception! To see this, first note that regular hexagons tessellate
the entire plane. But $\sin_{3}$ is already defined on one hexagon $\Omega
_{3}$, each side of which is mapped (by the boundary extension of $\sin_{3}$)
into a line. Schwarz reflection principle then extends $\sin_{3}$ analytically
over $%
\mathbb{C}
$ except at exactly half of the vertices, the $\omega_{3}^{k}P_{3}$'s among
them.\footnote{At the other half of the vertices, i.e., the translates of the
$\omega_{3}^{k}A_{3}$'s, the extended $\sin_{3}$ is bounded and hence analytic
(by Riemann's principle of removable singularity.)} These exceptional vertices
are \textit{isolated} singularities toward which $|\sin_{3}|\rightarrow\infty$
and hence they are the poles; among them, the $\omega_{3}^{k}P_{3}$'s are
nearest the origin.

Note that $R_{3}=\pi_{3}/2$ (the length of the interval on which the real
$\sin_{3}$ is originally defined) but $R_{n}<\pi_{n}/2$ (unexpectedly)\ for
$n>3$.

We leave a few questions for the reader to ponder.

\begin{exercise}
\label{Periodicity} Scrutinize the action of $\sin_{3}:%
\mathbb{C}
\rightarrow%
\mathbb{C}
\cup\{\infty\}$ to see that it is doubly periodic with periods $3\pi_{3}/2$
and $\left(  3\pi_{3}/2\right)  \omega_{3}$. (In fact, $\sin_{3}$ is known as
a Dixon's elliptic function; see \cite{Dixon}.)
\end{exercise}

\begin{exercise}
As squares also tessellate the plane and $\Omega_{4}$ is a square, one can ask
whether $\sin_{4}$ admits a meromorphic extension over $%
\mathbb{C}
$ in the same way $\sin_{3}$ does. (Answer: No! Why not?) What about
$(\sin_{4})^{2}$? (Answer:\ Yes! And it is also doubly periodic. This fact is
obtained in \cite{Levin} via a far less elementary approach.)
\end{exercise}

\subsection{Evaluating real integrals}

Many integrals arise from considering generalized trigonometric functions;
see, e.g., \cite{Poodiack}. Here, we consider two that yield to our analyses.

We saw (see Figure\thinspace\ref{Fn}) that $F_{n}[J]=[A_{n},P_{n})$ and
therefore%
\[
|A_{n}P_{n}|=\int_{1}^{\infty}\frac{1}{(t^{n}-1)^{(n-1)/n}}dt\text{\thinspace
;}%
\]
see (\ref{Integral}). As $[A_{n},P_{n})$ is a side of $\triangle OA_{n}P_{n}$
(see Figure\thinspace\ref{PI-n}), we have, for $n\geq3$,%
\[
\int_{1}^{\infty}\frac{1}{(t^{n}-1)^{(n-1)/n}}dt=\frac{\pi_{n}}{4}\sec
\frac{\pi}{n}\,\text{.}%
\]

\begin{exercise}
\label{Integral-OP}Show that $F_{n}[e^{i\pi/n}(I\cup J)]=[O,P_{n})$ and
therefrom deduce that%
\[
\int_{0}^{\infty}\frac{1}{(1+t^{n})^{(n-1)/n}}dt=\frac{\pi_{n}}{4}\sec
\frac{\pi}{n}\,\text{.}%
\]

\end{exercise}

\begin{exercise}
By the last exercise, $\sin_{n}[O,P_{n})$ is a ray. Apply Schwarz reflection
principle to relate the action of $\sin_{n}$ on the two triangles $\triangle
OA_{n}P_{n}$ and $\triangle OB_{n}P_{n}$. Further deduce that, for $z\in
\Omega_{n}$, $\sin_{n}(\omega_{n}z)=\omega_{n}\sin_{n}z$. (Hint: for the last
assertion, recall that the composition of two reflections across two
intersecting lines is a rotation.)
\end{exercise}

\subsection{Trisection of area between $C_{3}$ and its asymptote}

In the introduction, we asserted that \textit{the area between }$C_{3}%
$\textit{\ and its asymptote is trisected by the coordinate axes.} We argue
this based on properties of $\sin_{3}$.

Recall that, for real $t\in(F_{3}(-\infty),\pi_{3}/2)$, $\sin_{3}t$ is the
$y$-coordinate of the point $P(t)\in C_{3}$ such that the sector $OAP(t)$ has
area $t/2$ ($A$ being the point $(1,0)$). See Figure\thinspace\ref{C3} for
visual aid. If $P(t)$ marches along $C_{3}$ in the southeast direction, then
$\sin_{3}t$ runs monotonically through negative reals and tends to $-\infty$.
Therefore, to find the area of the unbounded\ sector that is the region in the
fourth quadrant between $C_{3}$ and its asymptote, it suffices to find the
first negative number $\alpha$ such that $\sin_{3}\alpha=-\infty$. A glance at
Figure\thinspace\ref{Domains(3,4,6)} tells us that $\alpha=\omega_{3}P_{3}$
and the area sought is%
\[
|\alpha|/2=|OP_{3}|/2=|OA_{3}|/2=\pi_{3}/4\text{\thinspace,}%
\]
which is exactly the area bounded by $C_{3}$ in the first quadrant.

Knowing that the area between $C_{3}$ and its asymptote is $3(\pi_{3}/4)$, we
can interpret the fact that $\sin_{3}(t+3\pi_{3}/2)=\sin_{3}t$ (found in
Exercise \ref{Periodicity}) by imagining that the two \textquotedblleft
ends\textquotedblright\ of $C_{3}$ meet at $\infty$, thereby turning $C_{3}$
into a closed curve, so that, when the sectorial area $t/2$ increases by
$3(\pi_{3}/4)$, $P(t)$ returns to the starting point.

\begin{exercise}
Give an alternative proof of the area trisection result using the integral
identity in Exercise \ref{Integral-OP}. (Hint:\ Find $\lim_{y\rightarrow
-\infty}F_{3}(y)$.)
\end{exercise}

\begin{exercise}
Does the area trisection claim remain true for $C_{2m+1}$ with $m>1$?
\end{exercise}

\section{Postscript}

We close with a few remarks.

The action of $F_{n}$ on $\overline{V}_{n}$ is reminiscent of that of a
Schwarz-Christoffel mapping, which maps the open upper-half plane $H$
conformally onto the interior of a convex polygon. We make the relation a bit
more exact. Define $\sqrt{V_{n}}$ to be $\{re^{i\theta}:r>0;\theta\in
(0,\pi/n)\}$. By Exercise \ref{Integral-OP}, $F_{n}$ maps $\sqrt{V_{n}}$
conformally onto $\mathring{\triangle}OA_{n}P_{n}$, the interior of $\triangle
OA_{n}P_{n}$. There is a \textit{unique} conformal map $S_{n}:H\rightarrow
\mathring{\triangle}OA_{n}P_{n}$ whose continuous extension maps $%
\mathbb{R}
\cup\{\infty\}$ onto $\partial(\triangle OA_{n}P_{n})$ in such a way that
$S_{n}(0)=O$, $S_{n}(1)=A_{n}$, and $S_{n}(\infty)=P_{n}$. As the function
$z\mapsto z^{n}$ maps $\sqrt{V_{n}}$ conformally onto $H$, we have, by
uniqueness, $F_{n}(z)=S_{n}(z^{n})$ for $z\in\sqrt{V_{n}}$. With $S_{n}$
expressed by a Schwarz-Christoffel formula, this functional identity yields
another integral expression for \textquotedblleft$\arcsin_{n}$%
\textquotedblright\ (i.e. $F_{n}$) on $\sqrt{V_{n}}$.

Turning to $\sin_{n}$, recall that Riemann mapping theorem guarantees the
existence of a unique conformal equivalence $\varphi_{n}:\Omega_{n}%
\rightarrow\Sigma_{n}$ such that $\varphi_{n}(0)=0$ and $\varphi_{n}^{\prime
}(0)>0$. Through our work, we have found it, i.e., $\varphi_{n}=\sin_{n}$.

Lastly, we note that the approach we have taken is readily applicable to a
class of so-called generalized sine functions with two parameters.
Specifically, for integer $n\geq2$ and positive real $p>1$, the real function
$\sin_{p,n}$ defined by the equation
\[
x=\int_{0}^{\sin_{p,n}(x)}(1-t^{n})^{-1/p}dt
\]
can be analytically continued in a similar fashion. For this generalization,
see \cite{Ding}.

\end{document}